\documentclass[12pt]{article}

\usepackage{amsmath}
\usepackage{amssymb}
\usepackage{amsfonts}

\begin{document}
\title{ The Sequence of Partial Quotients of \\Continued Fractions Expansion of Any Real \\Algebraic Number of Degree 3 is Bounded}

\author{Jinxiang Li\\
\small School of Sciences,\\ \small Guangxi University for Nationalities,\\
\small Nanning 530006, P.R.China}

\date{\small E-mail: lijinxiang-270@gxun.edu.cn }

\maketitle

\abstract{ Let $\,\alpha$ be a real algebraic number of degree 3. In this paper, by using the cubic formula of the cubic equation and the properties of the continuous function$\,f(x)=(1+p\theta)^{x}$ in the$\,p-$adic number field$\,Q_{p},$ I prove that if rational fraction$\;p_{1}/q_{1}$ such that$\;|\alpha-p_{1}/q_{1}|<q_{1}^{-2-\tau} \,(\tau>0),$ then $\,q_{1}^{\tau}<C.$ (where $\;C=C(\alpha)$ is an effectively computable constant.) In particular, the sequence of partial quotients of continued fractions expansion of any real algebraic number of degree 3 is bounded.

 In fact, we can conjecture that the sequence of partial quotients of continued fractions expansion of any real algebraic number is bounded.

\textbf{Keywords:} rational approximation, continued fraction, consecutive convergents, greatest prime factor,$\,p-$adic number,$\,p-$adic number field$\,Q_{p},$$\,p-$adic valuation.

\textbf{MSC(2000):} O156. }

\vspace{0.5cm}

\textbf{1. Introduction:}\\

Let$\;\alpha$ is a real algebraic number of degree$\;n\geq 2,$ there is a computable number$\;c=c(\alpha)$ such that\vspace{-0.5cm}
\begin{eqnarray*}
|\alpha-p/q|>cq^{-n}.
\end{eqnarray*}
for all rational numbers$\;p/q .$ This follows directly from the definition of an algebraic number,as was shown by Liouville in 1843; Axel Thue[3] was the first to prove a stronger result when$\;n\geq 3;$ he showed that the inequality\vspace{-0.2cm}
\begin{eqnarray*}
|\alpha-p/q|<q^{-0.5n-1-\tau},\;\tau>0,
\end{eqnarray*}
has at most finitely many solutions$\;(p,q),(p,q)=1.$

A further improvment was made by Siegel[4] in 1921; he proved that \vspace{-0.2cm}
\begin{eqnarray*}
|\alpha-p/q|<q^{-n(s+1)^{-1}-s-\tau},\;\tau>0,\;1\leq s <n.
\end{eqnarray*}
has at most finitely many solutions$\;(p,q),(p,q)=1.$

A further weakening was made by Dyson[5] and Gelfond[6] independently in 1948. They proved that\vspace{-0.4cm}
\begin{eqnarray*}
|\alpha-p/q|<q^{-\sqrt{2n}-\tau},\;\tau>0.
\end{eqnarray*}
has at most finitely many solutions$\;(p,q),(p,q)=1.$

Finally in 1955 Roth[7] obtained the best result, he proved that\vspace{-0.2cm}
\begin{eqnarray*}
|\alpha-p/q|<q^{-2-\tau},\;\tau>0.
\end{eqnarray*}
has at most finitely many solutions$\;(p,q),(p,q)=1.$

We know that there are infinitely many $\;p/q,(p,q)=1$ with \vspace{-0.2cm}
\begin{eqnarray*}
|\alpha-p/q|<q^{-2}.
\end{eqnarray*}
For any given $\;\alpha,$ with degree $\,\deg \alpha\geq 3,$ It is still unknown whether is badly approximable, i.e. whether there exists a $\,c>0$ so that \vspace{-0.4cm}
\begin{eqnarray*}
|\alpha-p/q|>cq^{-2},
\end{eqnarray*}
for every rational$\,p/q.$ The conjecture[1] is that this holds for no algebraic$\,\alpha$ of degree$\,\geq 3.$

Another conjecture[1] is that the inequality \vspace{-0.2cm}
\begin{eqnarray*}
|\alpha-p/q|<1/q^{-2}(\log q)^{k},
\end{eqnarray*}
has only finitely many solutions$\,p/q$ for$\,k>1.$

In this paper, I prove the following

\textbf{Theorem :} Let$\;\alpha$ is a real algebraic number of degree$\;n=3,$ if the inequality \vspace{-0.2cm}
\begin{eqnarray}\label{L7}
|\alpha-p/q|<q^{-2-\tau},\;\;q,\,\tau>0,\,(p,q)=1.
\end{eqnarray}
has rational number solutions$\;p, q,$ then$\;q^{\tau}<C=C(\alpha)$ (where$\;C$ is an effectively computable constant).In particular, the sequence of partial quotients of continued fractions expansion of any real algebraic number of degree 3 is bounded.

The second part of the theorem is true because a property of continued fraction.i.e.,
 If $\,\left|p/q-\alpha\right|<1/2q^{-2},$ then $\,p/q$ is a convergent.

So we only need to prove the first part of the theorem.\\

\textbf{2. Preliminaries:}\\

   In this section,$\,\alpha$ is an arbitrary real number rather than a cubic algebraic number. We first give some basic properties of continuous fractions of real number$\,\alpha.$

   Let$\,p_{1}/q_{1},\;p_{2}/q_{2}\,(q_{1}<q_{2})$ are two consecutive convergents to$\,\alpha,$ since the convergents are alternately less and greater than$\,\alpha,$ we have\vspace{-0.2cm}
\begin{eqnarray}\label{L10}
\left|\frac{p_{1}}{q_{1}}-\alpha\right|+\left|\frac{p_{2}}{q_{2}}-\alpha\right|=\left|\frac{p_{1}}{q_{1}}-\frac{p_{2}}{q_{2}}\right|=\frac{1}{q_{1}q_{2}},\;p_{2}q_{1}-p_{1}q_{2}=\pm1.
\end{eqnarray}
write\vspace{-0.4cm}
\begin{eqnarray}\label{L20}
\varepsilon_{1}=\frac{p_{1}}{q_{1}}-\alpha=\frac{\pm 1}{q_{1}(wq_{1}+q_{0})}=\frac{\pm 1}{q_{1}^{2+\tau}},\;\;\varepsilon_{2}=\frac{p_{2}}{q_{2}}-\alpha=\frac{\mp 1}{q_{2}^{2+\sigma}},\;(\tau>0,\;\sigma>0)
\end{eqnarray}
Note that$\;w=[w]+w',\;0<w'<1,$ so$\;q_{2}=[w]q_{1}+q_{0}=wq_{1}+q_{0}-w'q_{1}=q_{1}^{1+\tau}-w'q_{1}=q_{1}^{1+\tau}\left(1-w'q_{1}^{-\tau}\right).$ hence\vspace{-0.4cm}
\begin{eqnarray}\label{L30}
q_{2}=q_{1}^{1+\tau}\left(1-w'q_{1}^{-\tau}\right).
\end{eqnarray}
We have by (\ref{L20}) and (\ref{L30})\vspace{-0.4cm}
\begin{eqnarray}\label{L40}
\varepsilon_{2}=-s\,\varepsilon_{1},\;\;s=q_{1}^{-\tau}q_{2}^{-\sigma}\left(1-w'q_{1}^{-\tau}\right)^{-2}.
\end{eqnarray}
It is clear that$\;s>0,$ for $\;0<w'<1,\,q_{2}^{\sigma}>1.$ so $\;s\rightarrow 0$ when$\;q_{1}^{\tau}\rightarrow \infty.$ we substitute (\ref{L20}) into (\ref{L10}) we have\vspace{-0.4cm}
\begin{eqnarray}\label{L50}
\frac{1}{q_{1}^{2+\tau}}+\frac{1}{q_{2}^{2+\sigma}}=\frac{1}{q_{1}q_{2}}
\end{eqnarray}
We substitute (\ref{L30}) into (\ref{L50}) we obtain\vspace{-0.3cm}
\begin{eqnarray}\label{L60}
\frac{1}{q_{1}^{2+\tau}}+\frac{1}{q_{2}^{\sigma}q_{1}^{2+2\tau}(1-w'q_{1}^{-\tau})^{2}}=\frac{1}{q_{1}^{2+\tau}(1-w'q_{1}^{-\tau})}
\end{eqnarray}
Multiplying two side of (\ref{L60}) by$\;q_{1}^{2+\tau}$ we obtain\vspace{-0.3cm}
\begin{eqnarray}\label{L70}
1+\frac{1}{q_{1}^{\tau}q_{2}^{\sigma}(1-w'q_{1}^{-\tau})^{2}}=\frac{1}{1-w'q_{1}^{-\tau}}
\end{eqnarray}
So (\ref{L70}) gives \vspace{-0.5cm}
\begin{eqnarray}\label{L80}
1+s=\frac{1}{1-w'q_{1}^{-\tau}},\;i.e.\;(1+s)(1-w'q_{1}^{-\tau})=1,
\end{eqnarray}
by (\ref{L40}). Combining (\ref{L30}) (\ref{L80}) we also have\vspace{-0.3cm}
\begin{eqnarray}\label{L90}
q_{2}q_{1}=q_{1}^{2+\tau}(1-w'q_{1}^{-\tau})=\pm\varepsilon_{1}^{-1}(1+s)^{-1}.
\end{eqnarray}

\textbf{Lemma 1[2]:} Let $\,f\in Z[x,y]$ be a binary form such that among the linear factors in the factorization of $\,f$ at least three are distinct. Let $\,d$ be a positive integers. and $\,P$ be the greatest prime factor of $\;f(x,y).$ Then for all pairs of integers $\,x,y $ with $\,(x,y)=d,$\vspace{-0.3cm}
\begin{eqnarray}\label{L100}
P >> log log X,
\end{eqnarray}
where $\,X=max(|x|,|y|)>e$ and the possible constant implied by the $\,>>$ symbol only depends on $\,f$ and $\,d$ and is effectively computable.

We need some properties of the $p$-adic numbers field$\;Q_{p} $ and  integral domain$\;Z_{p}[\sqrt{-3}]. $ 

\textbf{Lemma 2[8]:} Let$\,\sqrt{-3}\not\in Z_{p}.$ For $\,\theta\in Z_{p}[\sqrt{-3}]$ the function$\,f(x)=(1+p\theta)^{x}$ is a continuous function from $\,Z_{p}$ to $\,1+pZ_{p}[\sqrt{-3}].$
i.e. For any $\,x\in Z_{p},$ there exists $\,\beta \in Z_{p}[\sqrt{-3}] $ such that\vspace{-0.3cm}
\begin{eqnarray}\label{L110}
(1+p\theta)^{x}=1+p\beta,\;\;\beta \in Z_{p}[\sqrt{-3}].
\end{eqnarray}
\\

\textbf{3. Proof of the theorem :}\\

Suppose$\;\alpha$ is real algebraic of degree 3 satisfies an equation\vspace{-0.2cm}
\begin{eqnarray}\label{L200}
ax^{3}+bx^{2}+cx+d=0,\;a\neq 0,\,a,\,b,\,c,\,d \in Z
\end{eqnarray}
and its discriminant on $\,x$ is\vspace{-0.4cm}
\begin{eqnarray}\label{L220}
\Delta=-27\,{a}^{2}{d}^{2}+18\,adcb+{b}^{2}{c}^{2}-4\,{b}^{3}d-4\,{c}^{3}a
\end{eqnarray}
Notice that if$\,q^{\tau}\geq 2,$ then $\,p/q$ is convergent to $\,\alpha$ in (\ref{L7}). So we just need to consider the convergent to $\,\alpha.$
Let $\,p_{1}/q_{1},\;p_{2}/q_{2}\,(q_{1}<q_{2})$ are two consecutive convergents to $\,\alpha,$ we have\vspace{-0.4cm}
\begin{eqnarray}\label{L240}
\alpha=\frac{p_{2}\beta+p_{1}}{q_{2}\beta+q_{1}}
\end{eqnarray}
Substitute this in (\ref{L200}) we obtain\vspace{-0.2cm}
\begin{eqnarray}\label{L260}
A\beta^{3}+B\beta^{2}+C\beta+D=0
\end{eqnarray}
where\vspace{-0.4cm}
\begin{eqnarray}\label{L280}
\begin{array}{rl}
A=&ap_{2}^{3}+bp_{2}^{2}q_{2}+cp_{2}q_{2}^{2}+dq_{2}^{3};\\
B=&\left(3\,ap_{2}^{2}+2\,bp_{2} q_{2}+c\,q_{2}^{2}\right)p_{1}+ \left( b\,p_{2}^{2}+2\,cp_{2}+3\,dq_{2}^{2}\,q_{2}\right)q_{1};\\
C=&\left(3\,ap_{1}^{2}+2\,bp_{1} q_{1}+c\,q_{1}^{2}\right)p_{2}+ \left( b\,p_{1}^{2}+2\,cp_{1}+3\,dq_{1}^{2}\,q_{1}\right)q_{2};\\
D=&a{{p_{1}}}^{3}+b{{p_{1}}}^{2}{q_{1}}+c{p_{1}}\,{{q_{1}}}^{2}+d{{q_{1}}}^{3}.
\end{array}
\end{eqnarray}
discriminant of (\ref{L260}) on$\;\beta $ is\vspace{-0.3cm}
\begin{eqnarray}\label{L300}
-27\,A^{2}D^{2}+18\,ABCD+B^{2}C^{2}-4\,B^{3}D-4\,AC^{3}=(p_{2}q_{1}-p_{1}q_{2})^{6}\Delta=\Delta.
\end{eqnarray}

We will prove that $\;q_{1}^{\tau}$ cannot be too large in (\ref{L300}).

We may write (\ref{L300}) as a cubic equation on$\;B$ ($\,\Delta$ fixed)\vspace{-0.2cm}
\begin{eqnarray}\label{L320}
-4\,DB^{3}+C^{2}B^{2}+18\,C A D B-27\,A^{2}D^{2}-4\,AC^{3}-\Delta=0.
\end{eqnarray}
Using the cubic formula of the cubic equation, we can decompose the left side of (\ref{L320}) into three factors $\,T_{1},\,T_{2}$ and $\,T_{3},$  therefore\vspace{-0.1cm}
\begin{eqnarray}\label{L410}
T_{1}\,T_{2}\,T_{3}=0,
\end{eqnarray}
where\vspace{-0.5cm}
\begin{eqnarray}\label{L420}
\begin{array}{rl}
T_{1}=&12DB-C^{2}-\sqrt[3]{E+12\sqrt{F}}-\sqrt[3]{E-12\sqrt{F}};  \\
T_{2}=&12DB-C^{2}-\omega\sqrt[3]{E+12\sqrt{F}}-\omega^{2}\sqrt[3]{E-12\sqrt{F}};\\
T_{3}=&12DB-C^{2}-\omega^{2}\sqrt[3]{E+12\sqrt{F}}-\omega\sqrt[3]{E-12\sqrt{F}}.
\end{array}
\end{eqnarray}
where$\;\omega\neq1,\;\omega^{3}=1,$ and
\begin{eqnarray}\label{L440}
\begin{array}{rl}
E=&C^6-108D^{2}\,(54 A^2 D^2+2\Delta +5 A C^3).\\
F=&3D^{2}\,(324D^{2}C^{6}A^{2}-8748D^{4}C^{3}A^{3}-4C^{9}A+78732D^{6}A^{4}\\
&+5832D^{4}A^{2}\Delta-C^{6}\Delta+108D^{2}\Delta^{2}+540D^{2}\Delta AC^{3}).
\end{array}
\end{eqnarray}

We shall prove that\\
\textbf{Lemma 3 :} If $\;q_{1}^{\tau}$ sufficiently large, then$\;T_{1}\neq 0,\;T_{2}\neq 0,\;T_{3}\neq 0.$ so (\ref{L410}) or  (\ref{L320}) is impossible.

It immediately follows that the theorem is true from the lemma 3.

   First, we have \vspace{-0.4cm}
\begin{eqnarray*}
p_{2}=q_{2}(\alpha+\varepsilon_{2}),\;p_{1}=q_{1}(\alpha+\varepsilon_{1}),\;\varepsilon_{2}=-s \,\varepsilon_{1}.
\end{eqnarray*}
by (\ref{L20}) and (\ref{L40}). Substitute this in (\ref{L280}) we obtain\vspace{-0.4cm}
\begin{eqnarray}\label{L460}
\begin{array}{rl}
A=& q_{2}^{3}\,u\,\varepsilon_{1}\,s\,(-1+\mu {s}\varepsilon_{1}-\delta {s}^{2}\varepsilon_{1}^{2});\\
B=& q_{2}^{2}q_{1}\,u\,\varepsilon_{1}\,(1-2\,s+\mu \left( -2+s \right) s\varepsilon_{1}+3\,\delta {s}^{2}\varepsilon_{1}^{2});\\
C=&q_{2}q_{1}^{2}\,u\,\varepsilon_{1}\,(2-s-\mu \left( -1+2\,s \right) \varepsilon_{1}-3\,\delta s\varepsilon_{1}^{2});\\
D=&q_{1}^{3}\,u\,\varepsilon_{1}\,(1+\mu \varepsilon_{1}+\delta \varepsilon_{1}^{2}).
\end{array}
\end{eqnarray}
where$\;u=3a\alpha^{2}+2b\alpha+c,\;v=3a\alpha+b,\;\mu=vu^{-1},\;\delta=au^{-1}.$

Note that$\,s\rightarrow 0$ also a fortiori $\,\varepsilon_{1}\rightarrow 0$ when$\;q_{1}^{\tau}$ sufficiently large by (\ref{L30}) and (\ref{L40}).

In order to prove that$\,T_{2}\neq 0,T_{3}\neq 0,$ we only need to prove $\,F>0$ when$\;q_{1}^{\tau}$ sufficiently large.

Combining (\ref{L300}), We can simplify the second equality of the (\ref{L440}) into
\begin{eqnarray}\label{L480}
F=3\,D^{2} \left( 48\,D^{2}{B}^
{2} -288\,ACD^{2}-{C}^{4}-8\,B{C}^{2}D\right)  \left( 9\,A C D-6\,{B}^{2}D+B{C}^{2} \right) ^{2}.
\end{eqnarray}

Substitute (\ref{L460}) in the first bracket of the (\ref{L480}), and note that $\,s>0,$ we obtain\vspace{-0.2cm}
\begin{eqnarray}\label{L500}
\begin{array}{rl}
&48\,D^{2}{B}^{2} -288\,ACD^{2}-{C}^{4}-8\,B{C}^{2}D\\
=&q_{2}^{4}q_{1}^{8}u^{4}\varepsilon^{4}(8^{3}\,s-192\,s^{2}+1280\,\mu\,s \,\varepsilon_{1}+ \cdots)>0.
\end{array}
\end{eqnarray}
So $\;F>0$ for sufficiently small$\;s,$ hence $\,\sqrt[3]{E+12\sqrt{F}},\,\sqrt[3]{E-12\sqrt{F}}$ are two unequal real numbers. Note that$\,\omega^{2}+\omega+1=0$, hence$\;T_{2}\neq 0,\;T_{3}\neq 0.$

Now let's prove $\;T_{1}\neq 0.$

Note that (\ref{L440}), the first equality of (\ref{L420}) can be rewritten as\vspace{-0.2cm}
\begin{eqnarray}\label{L520}
T_{1}=&12DB-C^{2}-C^{2}\sqrt[3]{1-12D\xi}-C^{2}\sqrt[3]{1-12D\bar{\xi}}
\end{eqnarray}
where\vspace{-0.4cm}
\begin{eqnarray}\label{L540}
\xi=C^{-6}(9DE'+\sqrt{3F'}),\;
\bar{\xi}=C^{-6}(9DE'-\sqrt{3F'}).
\end{eqnarray}
and \vspace{-0.4cm}
\begin{eqnarray}\label{L560}
\begin{array}{rl}
E'=&54 A^2 D^2+2\Delta+5 A C^3,\\
F'=&\left( 48\,D^{2}{B}^{2}-288\,ACD^{2}-{C}^{4}-8\,B{C}^{2}D\right)\left(9\,A C D-6\,{B}^{2}D+B{C}^{2}\right)^{2}.
\end{array}
\end{eqnarray}
Let$\;p>|3\Delta|$ be the greatest prime factor of$\;D,$ by the forth equality of (\ref{L280}) and lemma 1 and for sufficiently large$\;q_{1}.$  Also$\;(p,C)=1,$ by (\ref{L300}). Otherwise$\;p|\Delta $ which is impossible.

For the above determined prime number $\,p,$ we apply lemma 2 to consider (\ref{L520}) in the $p$-adic numbers field$\;Q_{p}.$

Note that$\,p|D,\,p\not|C,$ $\,1/2\in Z_{p},\;1+8DC^{-4}N \in 1+pZ_{p},$ so $\, M C^{2}\sqrt{1+8DNC^{-4}} \in Z_{p}.$\vspace{-0.2cm}
\begin{eqnarray}\label{L580}
\sqrt{3F'}=&M C^{2}\sqrt{1+8DNC^{-4}}\cdot\sqrt{-3} \in Z_{p}[\sqrt{-3}] .
\end{eqnarray}
where$\,M=\left|9\,A C D-6\,{B}^{2}D+B{C}^{2}\right|,\; N=B C^{2}+36\,A C D-6\,B^{2}D.$

Whether$\,\sqrt{-3}\in Z_{p} $ or$\,\sqrt{-3}\not\in Z_{p},$ $\,\sqrt{3F'}\in Z_{p}[\sqrt{-3}],$ and $\,v_{p}(\xi)\geq 0,\,v_{p}(\bar{\xi})\geq 0$ from(\ref{L540}). Therefore $\,12D\xi\in pZ_{p}[\sqrt{-3}],\,12D\bar{\xi}\in pZ_{p}[\sqrt{-3}].$ Note that$\,1/3\in Z_{p},$ hence there exist $\,\beta,\bar{\beta} \in Z_{p}[\sqrt{-3}] $ such that\vspace{-0.2cm}
\begin{eqnarray}\label{L940}
\sqrt[3]{1-12D\xi}=1+p\beta,\; \sqrt[3]{1-12D\bar{\xi}}=1+p \bar{\beta}.
\end{eqnarray}
by the lemma 2. Therefore we get
\begin{eqnarray}\label{L960}
C^{2}+C^{2}\sqrt[3]{1-12D\xi}+C^{2}\sqrt[3]{1-12D\bar{\xi}}= 3\,C^{2}+p\,C^{2}\,(\beta+\bar{\beta}).
\end{eqnarray}
Since $\;p\not|3C^{2},$ (\ref{L960}) gives$\,v_{p}\left(C^{2}+C^{2}\sqrt[3]{1-12D\xi}+C^{2}\sqrt[3]{1-12D\bar{\xi}}\right)=0.$
Note that$\,v_{p}(12DB)\geq 1,$ we have\vspace{-0.2cm}
\begin{eqnarray}\label{L980}
\begin{array}{rl}
v_{p}(T_{1})=&v_{p}\left(12DB-C^{2}-C^{2}\sqrt[3]{1-12D\xi}-C^{2}\sqrt[3]{1-12D\bar{\xi}}\right)\\
            =&v_{p}\left(12DB-3C^{2}-p\,C^{2}\,(\beta+\bar{\beta})\right)\\
            =& min\left\{v_{p}(12DB),v_{p}(3C^{2}+p\,C^{2}\,(\beta+\bar{\beta})\right\}\\
            =&v_{p}\left(3C^{2}+p\,C^{2}\,(\beta+\bar{\beta})\right)=0.
\end{array}
\end{eqnarray}
Hence $\;T_{1}\neq 0.$

From the above discussion, we can see that the upper bound of $\,q_{1}^{\tau}$ can be calculated effectively.

This completes the proof of the lemma 3.\\

\textbf{Acknowledgments:} I would like to thank Dr.Levent Alpoge for his careful reading of the original draft(arXiv:1904.09392[pdf.ps.other]math.NT) of this article, examining every step of calculations, and numerous helpful suggestions.
 \\

\textbf{References }\\
1. Schmidt W.M., Diophantine approximations and diophantine equations,\\
   Lecture Notes in Math.1467,Springer-Verlag,Berlin,Heidelberg,New York,1991. \\
2. A.Baker and D.W.Masser, Transcendence Theory:Advances and Applications.(1977).\\Academic Press London New York San Francisco.\\
3. A.Thue. \"{U}ber Ann\"{a}herungswerte algebraische Zahlen.\\J.reine angew.Math.,135(1909),284-305.\\
4. C.L.Siegel. Approximation algebraische Zahlen.Math.Z.,10(1921),173-213.\\
5. J.Dyson. The approximations to algebraic numbers by rationals.\\Acta Math.79(1947),225-240.\\
6. O.Gelfond. Transcendental and Algebraic Nubers.(Russian).\\English transl.(1969),Dover Publications,New York.\\
7. K.F.Roth. Rational approximations to algebraic numbers.\\Mathematika,2(1955),1-20.Also Corrigendum.Mathematika,2(1955),168.\\
8. K.Mahler. Introduction to p-adic Numbers and Their Functions,Cambridge University Press,1973.
\end{document}